\documentclass[a4paper,11pt,leqno]{amsart}

\usepackage{amssymb,latexsym}
\usepackage{graphicx}
\usepackage[all]{xy}
\usepackage{psfrag}

\usepackage[T1]{fontenc}
\usepackage[latin1]{inputenc}

\newcommand{\defop}[1]{\expandafter\def\csname #1\endcsname
{\mathsf{#1}}}

\newcommand{\C}{{\mathcal C}}
\newcommand{\E}{{\mathcal E}}

\newcommand{\eps}{\varepsilon}

\newcommand{\toto}{\longrightarrow}
\newcommand{\labelarrow}[2]{\buildrel{\scriptstyle#2}\over{#1}} %label#2/flèche #1
\newcommand{\iso}{\labelarrow{\toto}{\sim}}

\newcommand{\rdual}[1]{{#1^{\vee}}}
\newcommand{\ldual}[1]{{{}^{\vee}#1}}

\newcommand{\id}{\mathop{\mathsf{1}\kern 0pt}\nolimits}
\defop{End}
\defop{Hom}
\defop{Aut}
\defop{tr}
\defop{Ob}
\defop{vect}
\defop{kar}
\defop{Fun}
\defop{ann}

\newcommand{\Zed}{{\mathbb Z}}

\newtheorem{THEO}{Theorem}
\newtheorem{PROP}{Proposition}
\newtheorem{COR}{Corollary}
\newtheorem{LEM}{Lemma}

\newenvironment{DEF}{\noindent{\sc Definition.}}{\smallskip}
\newenvironment{NOT}{\noindent{\sc Notation.}}{\smallskip}
\newenvironment{REM}{\noindent{\sc Remark.}}{\smallskip}

\newenvironment{PROOF}{\noindent{\sc Proof.}}{$\,\Box$\smallskip}
\newenvironment{EX}{\noindent{\sc Example.}}{}

{\end{verse}\end{flushleft}}

\newcommand{\draw}[1]{\raisebox{-.45\height}{\includegraphics[scale=1]{#1.eps}}}
\newcommand{\minidraw}[1]{\raisebox{-.45\height}{\includegraphics[scale=.75]{#1.eps}}}
\newcommand{\Draw}[1]{\raisebox{-.45\height}{\includegraphics[scale=1.5]{#1.eps}}}

\defop{Im}
\defop{Tr}
\defop{Sph}
\defop{sl}
\defop{RPB}
\defop{PB}
\defop{RStL}
\defop{Coalg}
\defop{Comod}
\defop{mod}
\defop{can}

\newcommand{\SL}{\mathop{\mathsf{SL}\kern 0pt}\nolimits}
\newcommand{\PGL}{\mathop{\mathsf{PGL}\kern 0pt}\nolimits}
\newcommand{\SU}{\mathop{\mathsf{SU}\kern 0pt}\nolimits}
\newcommand{\rep}{\mathop{\mathsf{rep}\kern 0pt}\nolimits}
\newcommand{\SOV}{\mathop{\mathsf{SOV}\kern 0pt}\nolimits}
\newcommand{\Sov}{\mathop{\mathsf{Sov}\kern 0pt}\nolimits}

\defop{Tang}
\defop{eTang}
\defop{Turb}
\defop{StL}
\defop{Diag}
\defop{Br}
\defop{comod}

\newcommand{\D}{{\mathcal D}}

\newcommand{\sst}[1]{{\scriptscriptstyle{#1}}}

\newcommand{\del}{\partial}

\newcommand{\m}{^{-1}}

\newcommand{\CHI}[1]{{
{\raise 2pt\hbox{$\chi$}}\mkern -3mu
{}_{#1}
}}

\newcommand{\rvec}{\overrightarrow}
\newcommand{\lvec}{\overleftarrow}
\newcommand{\cl}[1]{\{{#1}\}}

\newcommand{\plus}{[+]}
\newcommand{\minus}{[-]}
\newcommand{\pt}{[0]}

\newcommand{\entrylabells}[1]{{(\bf #1)}}
\newenvironment{axioms}[1]%
 {\begin{list}{}{\settowidth{\labelwidth}{\entrylabells{#1}}%
  \settowidth{\labelsep}{\ \,}
  \setlength{\topsep}{0pt}
  \setlength{\itemsep}{0pt}
  \setlength{\leftmargin}{\labelwidth}\addtolength{\leftmargin}{\labelsep}%
  }}{\end{list}}

\newtheorem*{PROP1}{Proposition 1}
\title{Double Braidings, Twists and Tangle Invariants}
\author{A.~Bruguières}

\begin{document}

%\psfrag{=}{$=$}

\maketitle

\begin{abstract} A tortile (or ribbon) category defines invariants of ribbon (framed)
links and tangles. We observe that these invariants, when
restricted to links, string links, and more general tangles which
we call \emph{turbans}, do not actually depend on the braiding of
the tortile category. Besides duality, the only pertinent data for
such tangles are the double braiding and twist. We introduce the
general notions of twine, which is meant to play the rôle of the
double braiding (in the absence of a braiding), and the
corresponding notion of twist. We show that the category of
(ribbon) pure braids is the free category with a twine (a twist).
We show that a category with duals and a self-dual twist defines
invariants of stringlinks. We introduce the notion of
\emph{turban} category, so that the category of turban tangles is
the free turban category. Lastly we give a few examples and a
tannaka dictionary for twines and twists.
\end{abstract}

\begin{flushleft}
\begin{verse}`Just the place for a Snark!', the Bellman cried,\\ As he
landed his crew with care;\\ Supporting each man at the top of the
tide\\ With a finger entwined in his hair. \end{verse}
\end{flushleft}
\hfill {\small Lewis Carrol, \emph{The Hunting of the Snark}}

\tableofcontents

\thispagestyle{empty}

%\newpage
\section*{Introduction}

It is now well understood that certain categorical notions are
very closely related to low dimensional topology. For instance,
braids form a braided monoidal category, and the category of
braids is the free braided category. The category $\Tang$ of
oriented ribbon tangles is a tortile (or ribbon) category
\cite{JS1}, and indeed, it has been proved by Shum \cite{Shum}
(also \cite{tu1}) that $\Tang$ is the free tortile category. This
theorem is a powerful tool for constructing invariants of ribbon
links in $S^3$, since ribbon links up to isotopy are the
endomorphisms of the unit object in $\Tang$. Via Kirby calculus,
Shum's theorem underlies the construction of the
Reshetikhin-Turaev invariants of closed $3$-manifolds. Kirby
calculus can also be used to describe cobordisms of $3$-manfolds
in terms of certain tangles, and this allowed Turaev to construct
a TQFT associated with a modular category \cite{tu1}.

%Recall that a (strict) monoidal category is a category $\C$
%equipped with a tensor product, i. e. an associative functor
%$\otimes: \C \times \C \to \C$, and a unit object, i. e. and
%object $I$ such that $ ? \otimes I = \id_\C = I \otimes \, ?$. A
%braiding is a functorial isomorphism
%$$R_{X,Y}: X \otimes Y \iso Y \otimes X, \quad\mbox{$X,Y$ in
%$\C$}$$ satisfying

\medskip

The present work explores certain consequences of the following
observation. Let $\C$ be a tortile category. Recall that $\C$ is a
braided category with duals, and a (self-dual) twist $\theta$.
Denoting $c_{X,Y}: X \otimes Y \iso Y \otimes X$ the braiding,
define the \emph{double braiding} by $D_{X,Y}= c_{Y,X}c_{X,Y}$.
Notice that $\theta$  satisfies certain axioms where $c$ appears
only in the form of its double $D$, and conversely, $\theta$
determines $D$. It turns out that many significant notions
apparently related to $c$ actually depend only on $D$ or $\theta$.
The $S$-matrix, and the subcategory of transparent objects
\cite{Bru3}, which play an important role in the construction of
invariants of $3$-manifolds, are defined purely in terms of the
double braiding $D$. More surprisingly, the invariants of ribbon
links defined by $\C$ via Sum's Theorem do not depend on the
actual braiding, but only on $D$ (see Proposition 1, and section
5); and this result extends to a much larger class of ribbon
tangles, namely those whose linking matrix is diagonal $\mod 2$.
Since these tangles play an important role here, we give them a
name: we call them \emph{turban} tangles.

\medskip
All this suggests that the double braiding and the twist deserve
to be studied for their own sake, and that the universal property
of the category of tangles, that is, Shum's Theorem, should have
an analogue for the category of turban tangles.

The first step is to axiomatize  the notion of double braiding. We
observe that a double braiding satisfies certain formal properties
(TW0)-(TW2). An operator $D$ satisfying these properties will be
called a \emph{twine}. An \emph{entwined} category is a monoidal
category with a twine. The category of pure braids is the free
entwined category (section 2, Theorem 1). We also introduce a
general notion of twist, in such a way that the category of ribbon
pure braids is the free category with twist (section 3, Theorem
2).

In section 4, we bring duality into the picture, and find out that
a category with duals and (self-dual) twist defines invariants of
ribbon string links (theorem 3).

The heart of the matter is to extend these constructions to the
largest possible subcategory of the category of tangles. The
natural candidate is the subcategory generated by the twist and
duality: this is precisely the category of turban tangles
(proposition 1). In section 5, we define a \emph{turban category}
to be a category with a twist and good duals (sovereign
structure), satisfying certain additional conditions. We show that
the category of turban tangles is the free turban category
(theorem 4).

Section 6 gives a few examples of twines, twists and turban
categories, as well as the tannaka dictionary for twines and
twists.

\medskip
The definition of a turban category proposed in this paper is
certainly not definitive, but I believe that this notion could
lead to new topological invariants, including $3$-manifold
invariants and related TQFT's. The land of twines and twists is
`full of crags and chasms', and exploring it sometimes feels like
snark hunting. For instance, the fact that if $c$, $c$' are
braidings, then $c'_{Y,X}c_{X,Y}$ is a twine, came to me as a real
surprise! While so far I have few examples of twines or twists,
there are many indications that the class of entwined category is
much larger that the class of braided category.

\medskip
I must mention the fact that, after completion of this work, I was
informed that M. Staic recently obtained similar results in
\cite{Sta1}, where he constructs representations for the pure
braid group, and invariants of knots using Hopf algebra techniques
which correspond to the same categorical notions. Still, I hope,
the overlap is not so complete that one must conclude `the Snark
\emph{was} a Boojum, you see!'.

\medskip
I wish to thank Alexis Virelizier for many enlightening
discussions.

\section{Conventions and notations}

\subsection{Monoidal categories} Unless otherwise specified, all
categories will be small and all monoidal categories will be
strict. We will use Penrose Graphical calculus, with the ascending
convention: diagrams are to be read from bottom to top, {\it e.
g.} given $X\labelarrow{\to}{f}Y\labelarrow{\to}{g}Z$, we
represent $gf$ as
$$\draw{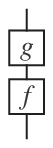}$$

If $\C$ is a monoidal category, with tensor product $\otimes$ and
unit object $I$, we denote $\otimes^n$ the $n$-uple tensor product
$$\displaylines {\C^n \toto \C\,,\cr (X_1, \dots X_n) \mapsto X_1
\otimes \dots \otimes X_n\,.}$$ In particular $\otimes^0 = I$,
$\otimes^1=\id_\C$ and $\otimes^2=\otimes$.

\medskip
Let $\C$ be a monoidal category. A \emph{duality of $\C$} is a
data $(X,Y,e,h)$, where $X$, $Y$ are objects, and $e: X \otimes Y
\to I$, $h: I \to Y \otimes X$ morphisms of $\C$, satisfying:
$$(e \otimes \id_X)(\id_X \otimes h)=\id_X \mbox{ and } (\id_Y \otimes e)(h\otimes \id_Y)=
\id_Y\,.$$ If $(X,Y,e,h)$ is a duality, we say that $(Y,e,h)$ is a
right dual of $X$, and $(X,e,h)$ is a left dual of $X$. If a right
or left dual of an object exists, it is unique up to unique
isomorphism.

\medskip
By \emph{monoidal category with right duals (resp. left duals,
resp. duals)}, we mean a monoidal category $\C$ where each object
$X$ admits a right dual (resp. a left dual, resp. both a right and
a left dual).

If $\C$ has right duals, we may pick a right dual
$(\rdual{X},e_X,h_X)$ for each object $X$ (the actual choice is
inocuous, in that a right dual is unique up to unique
isomorphism). This defines a monoidal functor
$$\rdual{?}: \C^o \to \C$$
where $\C^o$ denotes the category with opposite composition and
tensor product.

Similarly a choice of left duals $(\ldual{X},\eps_X,\eta_X)$ for
all $X \in \Ob\C$ defines a monoidal functor $\ldual{?}:\C^o \to
\C$.

\medskip
A \emph{(strict) sovereign structure on $\C$} is the choice, for
each object $X$, of a right dual $(X^*, e_X,h_X)$ and a left dual
$(X^*,\eps_X,\eta_X)$, with same underlying object $X^*$, in such
a way that $\ldual{?}=\rdual{?}$ as monoidal functors.
Essentially, left duals and right duals coincide.  By
\emph{sovereign category,} we mean a monoidal category with a
sovereign structure. This is an appropriate categorical setting
for a good notion of trace; however one must distinguish a left-
and a right trace $\tr_l$ and $\tr_r$. If $X$ is an object of $\C$
and $f \in \End(X)$,
$$\tr_l(f)=\eps_X(\id_{X^*}\otimes f)h_X,\quad \tr_r(f)=e_X(f\otimes \id_{X^*})\eta_X \quad\mbox{in $\End(I)$}\,.$$

\begin{DEF} Let $\C$ be a braided category, with braiding $c$; the
\emph{double braiding} is the functorial isomorphism
$$D_{X,Y}=c_{Y,X}c_{X,Y} : X \otimes Y \iso X \otimes Y\,.$$
\end{DEF}

A \emph{tortile category} is a monoidal braided category with
duals, equipped with a \emph{twist,} that is, a functorial
isomorphism $\theta_X: X \iso X$ ($X \in \Ob\C$) such that
$\theta_X \otimes \theta_Y = \theta_{X\otimes Y}D_{X,Y}$ and
$\theta_I=\id_I$. Moreover the twist is assumed to be
\emph{self-dual,} {\it i. e.}
$\theta_{\rdual{X}}=\rdual{\theta_X}$.

\medskip
If $\C$ is a tortile category, and if one makes the (inocuous)
choice of right duals $(\rdual{X}=X^*,e_X,h_X)$, there is a
canonical choice of left duals $(X^*,\eps_X,\eta_X)$ which defines
a sovereign structure. The self-duality of the twist implies that
the left- and right trace coincide (a property often referred to
as \emph{sphericity}).

\subsection{Tangles} We will often represent tangles by tangle
diagrams, which we view as drawings made up of the following
pictograms:
$$\draw{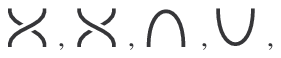}$$
called positive crossing, negative crossing, local max an local
min respectively (linked up by smooth arcs without horizontal
tangents).

Two tangle diagrams represent the same isotopy class of ribbon
tangles (also called framed tangles) if and only if one may be
obtained from the other by deformation and a finite number of
ribbon Reidemeister moves:
$$\displaylines{\mbox{(R0)} \;\draw{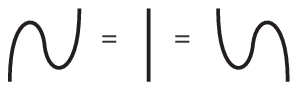}\;\,,\quad\mbox{(R1)}
\draw{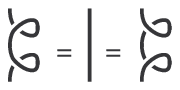}\,, \cr\quad \mbox{(R2)} \draw{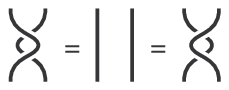}\,,
\draw{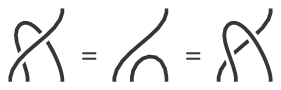}\,,\quad \mbox{(R3)} \draw{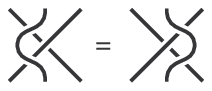}\,.}$$ Note that isotopy
of non-ribbon tangles is obtained by adding the Reidemeister move
$\raisebox{1mm}{\draw{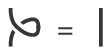}}\,$ to this list.

We will denote  $\cl{D}$ (and sometimes just $D$) the isotopy
class of ribbon tangles represented by a tangle diagram $D$.

\medskip
Let $D$ be a tangle diagram, $C$ a component of $D$.  We denote
$\Delta_C D$ the tangle diagram obtained from $D$ replacing $C$ by
$2$ parallel copies of $D$.

A tangle may be oriented, and/or coloured by elements of a set.

\medskip
We denote $\Tang$ the category of isotopy classes of oriented
ribbon tangles. This is a tortile category, whose objects are
words on the letters $\plus$ and $\minus$. We denote
$\Tang[\Lambda]$ the category of isotopy classes of oriented
ribbon tangles coloured by elements of the set $\Lambda$, which is
another tortile category. In $\Tang[\Lambda]$, we denote
$\plus_\lambda$ (resp. $\minus_\lambda$) the object $\plus$ (resp.
$\minus$) coloured by the element $\lambda \in \Lambda$. In the
oriented case, the point will be denoted $\pt$.

\subsection{Ribbon Tangles and Turbans Tangles}

Recall that the braiding $c$ and the twist $\theta$ of  the
category of ribbon tangles are defined by
$$c_{m,p}=\begin{matrix}{\draw{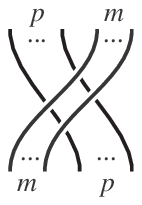}}, \quad \theta_n=\draw{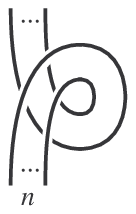}\,.,\end{matrix}$$
evaluation and coevaluation morphisms $e_n$ and $h_n$ being given
by $$e_n=\draw{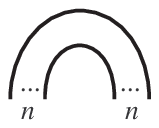}\,,\quad h_n=\draw{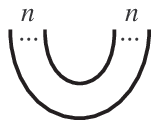}\,;$$ the ribbon
structure on $\Tang[\Lambda]$ is defined by the same tangles, with
appropriate orientation and $\Lambda$-colouring.

\begin{DEF} A ribbon tangle $T$ is \emph{turban} (resp. \emph{even}) if its linking
matrix is diagonal $\mod 2$ (resp. zero  $\mod 2$).
\end{DEF}

Turban tangles (resp. even tangles) form a monoidal subcategory of
$\Tang$ which we denote $\Turb$ (resp. $\eTang$).

The following proposition is the main motivation for the rest of
this work.

\begin{PROP}
The category $\Turb$ (resp. $\eTang$) is the smallest monoidal
subcategory of $\Tang$ having the same objects as $\Tang$, and
containing all evaluations, coevaluations and twists (resp. double
braidings).
\end{PROP}

We will prove proposition 1 in section 5.

\section{Twines and pure links}

\begin{DEF}
Let $\C$ be a monoidal category.  A  \emph{twine of $\C$} is an
automorphism $D$ of $\otimes$, that is, a functorial isomorphism
$$D_{X,Y}: X \otimes Y \iso X \otimes Y \qquad\mbox{($X, Y \in
\C$)}$$ satisfying the following axioms:
\begin{equation}\tag{DB0}D_{I,I}=\id_I\,;\end{equation}
\begin{equation}\tag{DB1}(D_{X,Y} \otimes \id_Z)D_{X
\otimes Y,Z}= (\id_X \otimes D_{Y,Z})D_{X,Y \otimes
Z}\,;\end{equation}
\begin{multline}\tag{DB2}(D_{X\otimes Y,Z} \otimes \id_T)(\id_X \otimes
D\m_{Y,Z}\otimes \id_T)(\id_X \otimes D_{Y,Z \otimes T})\\
= \quad(\id_X \otimes D_{Y,Z \otimes T})(\id_X \otimes
D\m_{Y,Z}\otimes \id_T)(D_{X\otimes Y,Z} \otimes
\id_T)\,.\end{multline} An \emph{entwined category} is a monoidal
category equipped with a twine. \end{DEF}

\medskip
If $\C$, $\C'$ are two entwined categories, with twines $D$, $D'$,
a strict entwined functor $F: \C \to \C'$ is a strict monoidal
functor $\C \to \C'$ such that for all $X, Y \in \Ob \C$,
$$F(D_{X,Y})=D'_{F X,F Y}\,.$$

\medskip
\begin{EX}
Let $\C$ be a monoidal category, and $c$ a braiding of $\C$. Then
$D_{X,Y}= c_{Y,X}c_{X,Y}$ (the double of $c$) is a twine of $\C$.
In particular, let $B$ be the category of braids, with its
canonical braiding $c$. Recall that the canonical braiding $c$ is
characterised by the fact that $c_{1,1}$ is the standard generator
of $B_2$. Let $D$ be the double of $c$. Let $\PB$ be the category
of pure braids. Then for any integers $m, p$, $D_{m,p}$ is a
morphism of $P$ and this defines a twine of $\PB$. We will
therefore consider $\PB$ as an entwined category.
\end{EX}

\medskip
\begin{REM} The examples admits of the following surprising generalization,
which was pointed out to me by A. Virelizier: if $c,c'$ are two
braidings in $\C$, then $D_{X,Y}=c'_{Y,X}c_{X,Y}$ is a twine.
\end{REM}

\medskip
Here are a few comments on the axioms.

The first two axioms (DB0) and (DB1) imply the following:

(a) $D_{X,I}=\id_X=D_{I,X}$.

(b) $(D\m_{X,Y} \otimes \id_Z) D_{X,Y \otimes Z}=D_{X \otimes
Y,Z}(\id_X \otimes D\m_{Y,Z})$ and $(\id_X \otimes D\m_{Y,Z}) D_{X
\otimes Y, Z}=D_{X, Y \otimes Z}( D\m_{X,Y} \otimes \id_Z)$.

It will be very convenient to depict $D_{X,Y}$, $D\m_{X,Y}$ as
follows:
$$D_{X,Y}= \draw{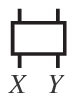}\, , \quad D\m_{X,Y}= \draw{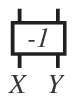}\,.$$
Similarly, let
$$\displaylines{D^f_{X,Y,Z}=\draw{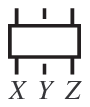}=(D\m_{X,Y} \otimes \id_Z) D_{X,Y \otimes Z}=D_{X\otimes Y,Z}(\id_X \otimes D\m_{Y,Z})
\,,\cr D^b_{X,Y,Z}=\draw{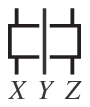}=(\id_X \otimes D\m_{Y,Z}) D_{X
\otimes Y, Z}=(\id_X \otimes D\m_{Y,Z})D_{X\otimes Y,Z}\,.}$$

Now (DB2) can be re-interpreted in a nice way. Indeed, composing
each side of (DB2) on the right by $(\id_X \otimes D\m_{Y,Z}
\otimes \id_T)$ and using (b), we obtain the sliding property:
$\draw{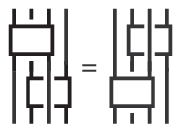}$.

\medskip
Notice that the notion of twine is invariant under left-right
symmetry (tensor product reversal) and under top-bottom symmetry
(composition reversal). In both cases front and back (i. e.
$\draw{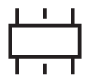}$ and $\draw{Dfxxx}$) are exchanged. In particular
central symmetry preserves front and back.
\medskip

The following theorem justifies, in a sense, the axioms for a
twine.

\begin{THEO} The category of pure braids is the universal entwined
category. More precisely, let $\C$ be an entwined category,
$\Lambda=\Ob\C$, and denote $\PB(\Lambda)$ the category of
$\Lambda$-coloured pure braids. There exists a unique strict
entwined functor  $\PB(\Lambda) \to \C$ sending $\pt_X$ ($X \in
\Ob\C$) to the object $X$ itself.
\end{THEO}

\subsection*{Proof of Theorem 1}

The proof relies on a presentation of the group of pure braids
$P_n$ by generators and relations, due to Markov \cite{Mar1}. (See
also \cite{Ver1}). Let $\sigma_i \in B_n$ ($1 \le i < n$) be the
standard generator:
$$\sigma_i= \draw{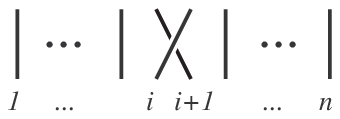}\,.$$
For $1 \le i < j \le n$, let $s_{i,j}=\sigma_{j-1} \dots \sigma_i
\sigma_i \dots \sigma_{j-1}$; pictorially:
$$s_{i,j}=\draw{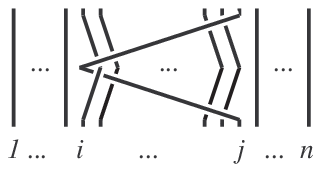}\,.$$ Then the $s_{i,j}$'s generate $P_n$,
subject to the Burau relations:

(Bu1) $s_{i,j}s_{k,l}=s_{k,l}s_{i,j}$ for $i < j< k< l$ or
$i<k<j<l$;

(Bu2)
$s_{i,j}s_{i,k}s_{j,k}=s_{i,k}s_{j,k}s_{i,j}=s_{j,k}s_{i,j}s_{i,k}$
for $i<j<k$;

(Bu3)
$s_{i,k}s_{j,k}s_{j,l}s\m_{j,k}=s_{j,k}s_{j,l}s\m_{j,k}s_{i,k}$
for $i<j<k<l$.

In the entwined category $\PB$ of pure braids, $s_{i,j} =
\id_{i-1} \otimes  D^f_{1,j-i-1,1}\otimes \id_{n-j}$.

\begin{PROP} Let $\C$ be an entwined category. There exists a unique group morphism
$$\displaylines{\PB_n \to \Aut(\otimes^n)\cr
P \mapsto [P]}$$ such that for all $X_1, \dots, X_n \in \Ob\C$ and
$1 \le i < j \le n$,
$$[\sigma_{i,j}]_{X_1, \dots, X_n} = \id_{X_1 \otimes \dots \otimes X_{i-1}}
\otimes \,D^f_{X_i, X_{i+1} \otimes \dots \otimes X_{j-1}, X_j}
\otimes \id_{X_{j+1} \otimes \dots \otimes X_n}\,.$$
\end{PROP}

\begin{PROOF} Since the $s_{i,j}$'s generate $P_n$, we only have
to check compatibility with the Burau relations.

Now the first case of (Bu1) is functoriality of the tensor
product, and the second case of (Bu2) is functoriality of
$D^f_{X,Y,Z}$ with respect to $Y$.

In order to check the other relations, we will have to perform
certain computations in an entwined category. Let us adopt the
notation:
$$\draw{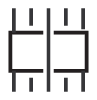}= \draw{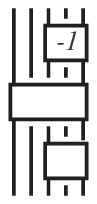}\,.$$

It is understood that each strand is coloured by an object of
$\C$, so this is an identity of morphisms of $\C$.

\begin{LEM}
The following identities hold in an entwined category:

a) $\draw{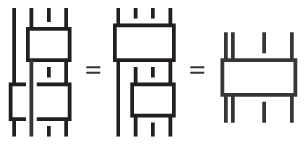}$; \quad b) $\draw{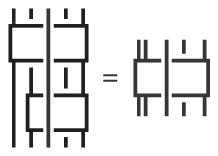}$; \quad c) $\draw{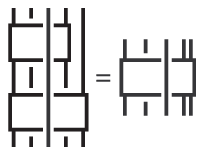}$;

d) $\draw{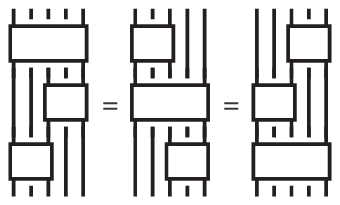}$; \quad e) $\draw{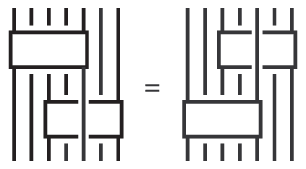}$.
\end{LEM}

\noindent {\it N. B.:} strings which are drawn very close
represent one entry coloured by the tensor product of the colours
of the strings.

\medskip
\begin{PROOF}
The computations would be very awkward in algebraic form; they are
much easier to conduct using Penrose graphical calculus. Here is a
sketch of the proof.

Assertions a) and b): the first identity of a) holds by
definition; the second results from the definition of
$\draw{Dfxxx}$ by straightforward computation, and implies b) by
definition of $\draw{FB}$.

Consider assertion e), and denote $X, A, Y, B, Z, C, T$ the
objects of $\C$ used to colour the seven strands, listed from left
to right. Then the case $A=B=C=I$ is just the sliding property,
which is a consequence of the twine axioms. Now using a), we
deduce e)  in the case $B=C=I$.

c) Using the definitions and elementary manipulations, assertion
c) can be easily reduced to assertion e) in the case $B=C=I$,
which we just proved, and the identity $\draw{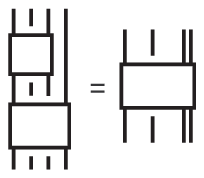}$, analogous to
the second identity of a).

Assertion e): the case $C=I$ can now be deduced from the case
$B=C=I$ using b). Hence the general case,  using a) and c).

Let us prove assertion d). By reason of symmetry, it is enough to
check the first identity. Now one computes easily
$$\draw{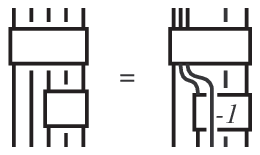}$$
and one concludes using e) and functoriality of the twine. Thus
ends the proof of the lemma. \end{PROOF}

Relations (Bu2) and (Bu3) are direct consequences of assertions d)
and e) of the lemma, hence the proposition.
\end{PROOF}

Now the lemma clearly defines a monoidal functor $\PB[\Lambda] \to
\C$ which sends $[X]$ to $X$ ($X \in \Ob\C$), the pure braid on
$n$ strands coloured by $X_1, \dots X_n$ to $[P]_{X_1, \dots,
X_n}$, and preserves the twine. Uniqueness results form the fact
that the $s_{i,j}$'s generate $\PB$.

\section{Twists and ribbon pure braids}

\begin{DEF} Let $\C$ be a monoidal category.
A \emph{twist of  $\C$} is an automorphism $\theta$ of $\id_\C$,
that is, a functorial isomorphism
$$\theta_X: X \iso X \qquad\mbox{($X \in \C$)}$$
satisfying the following axioms:
\begin{equation}\tag{TW0}\theta_{I}=\id_I\,;\end{equation}
\begin{multline}\tag{TW1}(\theta\m_{X\otimes Y} \otimes
\id_{Z\otimes T}) (\theta_{X \otimes Y \otimes Z} \otimes
\id_T)(\id_X \otimes \theta\m_{Y \otimes Z} \otimes \id_T)\\ (\id_X
\otimes \theta_{Y \otimes Z \otimes T})
(\id_{X\otimes Y} \otimes
\theta\m_{Z \otimes T})
 = (\id_{X\otimes Y} \otimes \,\theta\m_{Z \otimes T})
(\id_X
\otimes \theta_{Y \otimes Z \otimes T})\\
(\id_X \otimes \theta\m_{Y \otimes Z} \otimes \id_T) (\theta_{X
\otimes Y \otimes Z} \otimes \id_T)(\theta\m_{X\otimes Y} \otimes
\id_{Z\otimes T}) \,.\end{multline}
Graphically, axiom (TW1) may be represented as
$$\minidraw{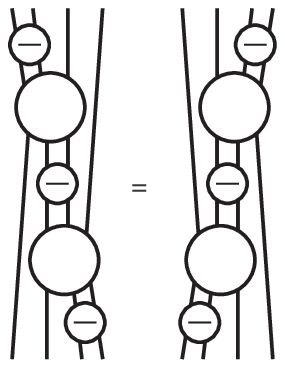}\,, \quad \mbox{with $\bigcirc=\theta$ and $\circleddash=\theta\m$.}$$
A \emph{twisted category} is a monoidal category equipped with a
twist.

If $\C$, $\C'$ are two twisted categories, with twists $\theta$,
$\theta'$, a strict twisted functor $F: \C \to \C'$ is a strict
monoidal functor satisfying for all $X \in \Ob\C$:
$$F(\theta_X)=\theta_{F X}\,.$$
\end{DEF}

\medskip
\begin{PROP} Let $\C$ be a monoidal category and $\theta$ an
automorphism of~$\id_\C$. Define an automorphism $D$ of $\otimes$
by
$$D_{X,Y}= (\theta\m_X \otimes \theta\m_Y)\theta_{X\otimes Y}\,.$$
Then $\theta$ is a twist if and only if $D$ is a twine.
\end{PROP}

\begin{PROOF}
By its very form, $D$ satisfies (DB1), and one checks easily that
(DB0) and (DB2) are equivalent respectively to (TW0) and (TW1).
\end{PROOF}

\medskip

As a result, a twisted category is canonically entwined, and a
strict twisted functor is entwined.

\begin{EX} Let $\C$ be a braided category, and let $\theta$ be a
balanced structure, that is an automorphism of $\id_\C$ satisfying
$$\theta_{X \otimes Y}= (\theta_X \otimes \theta_Y)
R_{Y,X}R_{X,Y}\,.$$ Then $\theta$ is a twist.

In particular, the category of ribbon braids is twisted, and so is
the category of ribbon pure braids. Moreover we have a canonical
group isomorphism $$(u,t_1, \dots, t_n) : \RPB_n \iso \PB_n \times
\Zed^n\,,$$ where $u$ denotes the forgetful morphism $\RPB_n \to
\PB_n$, and $t_i$ the self-linking number of the $i$-th component.

\medskip
\begin{REM} Let $\C$ be an entwined category, with twine $D$. Just like in the braided case,
there is a canonical way of adjoining a twist to $\C$. Indeed,
define a category $\tilde{\C}$ as follows. The objects of
$\tilde{\C}$ are data $(X,t)$, with $X \in \Ob\C$ and $t \in
\Aut(X)$. Morphisms from $(X,t)$ to $(X',t')$ are morphisms $f : X
\to X'$ in $\C$ such that $t'f=ft$. Define a tensor product on
$\tilde{\C}$, on objects, by
$$(X,t) \otimes (X',t')= (X \otimes X', (t \otimes t')D_{X,X'})\,,$$
and on morphisms, by the tensor product of $\C$.  One checks
easily that this makes $\tilde{\C}$ a strict monoidal category
(using axioms TW0 and TW1), and that setting $\theta_{(X,t)}=t$
defines a twist $\theta$ on $\tilde{\C}$ (using TW2). The
forgetful functor $\tilde{\C} \to \C$ is entwined, and this
construction is universal.
\end{REM}

\end{EX}

\begin{THEO}
Let $\C$ be a twisted category, $\Lambda=\Ob\C$, and
$\RPB(\Lambda)$ be the category of $\Lambda$-coloured ribbon pure
braids. There exists a unique strict twisted functor $$[?]
:\RPB[\Lambda] \to \C$$ sending $\pt_X$ ($X \in \Ob\C$) to the
object $X$ itself.
\end{THEO}

\begin{PROOF}
Clearly the image of the coloured ribbon pure braid $P[X_1, \dots,
X_n]$ can be no other than $(\theta^{t_1(P)}_{X_1} \otimes \dots
\theta^{t_n(P)}_{X_n}) [u(P)]_{X_1, \dots, X_n}$, and this defines
indeed a strict twisted functor.
\end{PROOF}

\begin{NOT} Let $\C$ be a twisted category; for $P \in \RPB_n$, we
let $[P] = (\theta^{t_1(P)} \otimes \dots \theta^{t_n(P)}) [u(P)]
\in \End(\otimes^n)$.
\end{NOT}

\section{Twists, duality, and invariants of string links}

\begin{DEF} Let $n$ be a non-negative integer. A \emph{(ribbon) string link}
on $n$ strands is an oriented (ribbon) tangle from $n\plus$ to
$n\plus$, without closed components, and such that the $i$-th
input is connected to the $i$-th output. \end{DEF}

\medskip
We denote $\RStL$ the monoidal subcategory of $\Tang$ whose
morphisms are isotopy classes of ribbon string links.

The category of ribbon pure braids $\RPB$ is naturally embedded as
a monoidal subcategory of $\RStL$.
\medskip

\begin{DEF}
Let $P$ be a ribbon string link on $n$ strands, and $1 < i < n$.
We define the \emph{$i$-th right contraction} of $P$ to be the
ribbon string link on $n-2$ strands $c_i(P)$ defined by
$${%\psfrag{P}{$P$}
c_i(P)= \begin{matrix}\overbrace{\hspace{8mm}}^{i-2}
\hspace{12mm}\overbrace{\hspace{8mm}}^{n-i-1}\cr\draw{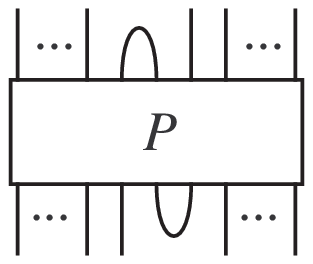}\cr
\underbrace{\hspace{8mm}}_{i-2}\hspace{12mm}\underbrace{\hspace{8mm}}_{n-i-1}\end{matrix}\,.}$$
\end{DEF}

We will now mimick this construction in a categorical setting,
using the notion of duality.

\begin{PROP}
Let $\C$ be a monoidal category, and let $\theta$ be a twist of
$\C$. Let $X$ be an object of $\C$ and $(Y,e,h)$ be a right dual
of $X$.

The following assertions are equivalent:

(i) $\theta_X$ and $\theta_Y$ are dual morphisms;

(ii) $\theta^2_X = (e D\m_{X,Y} \otimes \id_X)(\id_X \otimes h)$;

(ii') $\theta^2_Y= (\id_Y \otimes e D\m_{Y,X})(h\otimes \id_Y)$;

Moreover if they hold for one right dual of $X$, they hold for
all.

\end{PROP}

\begin{DEF} Let $\C$ be a monoidal category with right duals.
Let $\theta$ be a twist of $\C$. We say that $\theta$ is
\emph{self-dual} if for any object $X$ of $\C$ the equivalent
assertions of the previous proposition hold. \end{DEF}

\medskip

Let $\C$ be a twisted with right duals and a self-dual twist.
Assume that right duals are chosen.

\begin{NOT} Let $\C$ be a monoidal category with right duals.
Let $X,Y,Z$ be objects of $\C$. For $f \in \End(X \otimes Y
\otimes \rdual{Y} \otimes Y \otimes Z)$, let $$c_{X,Y,Z}(f)=
(\id_X \otimes \, e \otimes \id_Y \otimes \id_Z)f(\id_X \otimes
\id_Y \otimes \,h \otimes \id_Z)\,.$$

Now let $\phi \in \End(\otimes^n)$ and $1 < i < n$. Define $c_i
\phi \in \End(\otimes^{n-2})$ by
$$(c_i\ \phi)_{X_1, \dots, X_{n-2}}=c_{X_1 \otimes \dots \otimes X_{i-2},X_{i-1}, X_i \otimes \dots X_{n-2}}
\phi_{X_1, \dots, X_{i-1},\rdual{X_{i-1}},X_{i-1}, \dots,
X_n}\,.$$ Pictorially, \vspace{-5mm}
$$\psfrag{P}[c][c]{{\Large $\phi$}}
(c_i \phi)_{X_1, \dots,
X_{n-2}}=\begin{matrix}\hspace{8mm}\sst{{X_1 \dots X_{i-2}}
\hspace{6mm}\scriptscriptstyle{X_{i-1} \dots
X_{n-2}}}\hspace{5mm}\cr\draw{cip}\cr \hspace{10mm}\sst{X_1\,
\dots\quad X_{i-1}} \hspace{5mm}\sst{X_i \,\dots\,
X_{n-2}}\hspace{6mm}\end{matrix}\,.$$ Notice that $c_i\phi$ is in
fact independent of the choice of a right dual for $X_{i-1}$.
\end{NOT}

\begin{THEO}
Let $\C$ be a monoidal category with right duals and a self-dual
twist. There exists a unique way of associating to each isotopy
class of ribbon string link $P \in \RStL_n$ a functorial
endomorphim $\rvec{P} \in \End(\otimes^n)$ in such a way that:
\begin{enumerate}
\renewcommand{\labelenumi}{{\rm (\roman{enumi})}}
\item $\rvec{P}= [P]$ for any ribbon pure braid $P$;

\item $c_i \rvec{P}=\rvec{c_i P}$ for any $P \in \RStL_n$ and $1 <
i < n$ such that the $i$-th component of $P$ is trivial.
\end{enumerate}
\end{THEO}

\begin{COR}
Let $\C$ be a monoidal category with right duals and a self-dual
twist. Let $\Lambda=\Ob\C$ be a set and denote $\RStL[\Lambda]
\subset \Tang[\Lambda]$ the category of $\Lambda$-coloured ribbon
string links. There exists a canonical twisted functor
$$\rvec{?}: \RStL[\Lambda] \toto \C$$
which sends a coloured ribbon string link $P[X_1, \dots, X_n]$ to
$\rvec{P}_{X_1, \dots, X_n}$.
\end{COR}

If $\C$ is a twisted category with left duals, one may (by
left-right symmetry) associate with any ribbon string link $P \in
\RStL_n$ an element $\lvec{P} \in \End(\otimes^n)$. If both right-
and left duals exist, it is not at all clear whether
$\lvec{P}=\rvec{P}$. This suggests the following definition.

\medskip
\begin{DEF} Let $\C$ be a monoidal category with left and right
duals, and $\theta$ a twist of $\C$. We say that $\theta$ is
\emph{ambidextrous} if it is self-dual, and we have
$$\forall P \in \RStL, \,\lvec{P}=\rvec{P}\,.$$
If such is the case, we set $[P]=\lvec{P}=\rvec{P}$.
\end{DEF}

\medskip
When the twist is ambidextrous, we have $[c_i P]=c_i[P]$ for any
$P \in \RStL_n$ and any $1<i<n$.
\medskip

\begin{REM}
Theorem 3, while it provides a means of constructing invariants of
ribbon string links, has a serious drawback : it is not a
universal property, because the category of ribbon string links
has no duals. The aim of section 5 will be to mend this matter.
\end{REM}

%By \emph{compatible with contractions}, we mean the following. For
%

\subsection*{Proof of Theorem 3}

If $P$ is a ribbon pure braid, $\rvec{P}=[P]$ is already
well-defined. The point is now to see that a string link can be
obtained from a pure braid by a sequence of `nice' contractions.
This will at least show that $\rvec{?}$ is unique, and suggest a
construction for it. We then must check the coherence of this
construction, {\it i. e.} its independence from the choices made.

The main trick we use consists in `pulling a max to the top line'.
Let $D$ be a tangle diagram with a local max $m$, with $p$
outputs. We may write $$D=\draw{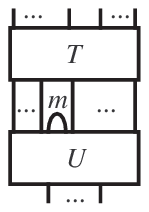}\,,$$ where $T$, $U$ are
tangle diagrams.

Let $i$ be an integer, $1 \le i \le n+1$. Let $j$ be the number of
strands to the left of $m$ on the same horizontal line, plus $1$.
Let $T'$ be a tangle diagram obtained from $T$ by inserting a new
component $C$ going from a point between the $(j-1)$-th and $j$-th
inputs of $T$ to a point between the $(i-1)$-th and $i$-th outputs
of $T$. We assume also that $C$ has no local extrema. Note that we
have $T=T'-C$. Let $T''=\Delta_C T'$ be the tangle diagram
obtained from $T'$ by doubling $C$. Set $$D'=\draw{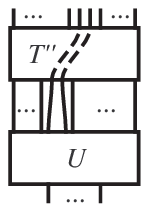}\,.$$

We say that \emph{$D'$ is obtained from $D$ by pulling $m$ to the
top in the $i$-th position (along the path $C$).}

One defines similarly the action of \emph{pulling a local min to
the bottom.}

Now let $D$ be a $n$-string link diagram, oriented from bottom to
top. We say that $D$ is \emph{right-handed} if all local extrema
point to the right.

Assume $D$ is right-handed. Pulling all local max to the top and
all local min to the bottom, one may obtain a pure braid diagram.
Here is an algorithm. Denote $m_i$ the number of local max (which
is equal to the number of local min) on the $i$-th component of
$D$. Let $m=m(D)=m_1 + \dots + m_n$ be the number of local max of
$D$. If $m(D)=0$, we are already done. Otherwise, chose $i$
minimal so that $m_i> 0$. Denote $c$ the $i$-th component, and let
$m$ be the first max, and $m'$ the first min you meet on $c$,
going from bottom to top. Pull $m$ to the top, in the $i$-th
position (just to the left of $c$), and $m'$ to the bottom, in the
$i+1$-th position (just to the right of $c$). Let $D'$ be the
diagram so constructed. Then $D'$ is a string link diagram, with
$m(D')=m(D)-1$. Moreover, $\cl{D}= c_{i+1} \cl{D'}$, and the
$(i+1)$-th component of $D'$ is unknotted. Repeated $m$ times,
this transformation yields a pure braid diagram $P$ with $n+ 2m$
threads, and we have
$$\cl{D}= c_{j_m} \dots c_{j_1} P\,,$$
where  $1 \le j_1 \le \dots \le j_m \le n$, and $j_k$ takes $m_i$
times the value $i+1$.

We therefore set $$\rvec{D}= c_{j_m} \dots c_{j_1} [P]\,,$$ and we
now proceed to show that this is independent of the choice made in
the construction of $P$, that is, the actual paths along which the
local extrema are being pulled.

\begin{LEM}
Let $P$, $P'$ be two pure braid diagrams:
$$P=\draw{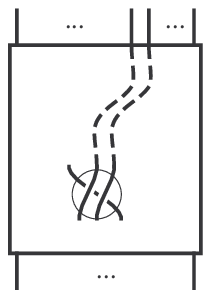} \,, \quad P'=\draw{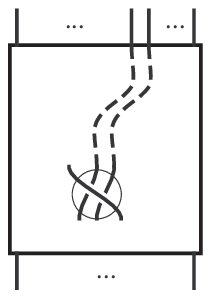}\,,$$
which differ only inside a circle. Inside the circle, the $i$-th
and $i+1$ strands pass respectively to the front and the back of
the $k$-th strand; above the circle, the $i$-th and $(i+1)$-th
strands run parallel. Let $\C$ be a entwined category, $X_1,
\dots, X_n$ objects of $\C$, and let $e : X_i \otimes X_{i+1} \to
I$ be any morphism. Let $E=\id_{X_1\otimes \dots \otimes X_{i-1}}
\otimes \, e \otimes \id_{X_{i+2}\otimes \dots \otimes X_n}$. Then
$E [P]_{X_1, \dots, X_n} = E [P']_{X_1, \dots, X_n}$.
\end{LEM}

\begin{PROOF}
We will use the following fact, which is an immediate consequence
of Proposition 1. If $A \in \PB_n$ and $1 \le i \le n$, construct
$\Delta_i P \in \PB_{n+1}$ by doubling the $i$-th strand of $P$.
Given $n+1$ objects $X_1, \dots, X_i , X'_i, X_{i+1}, \dots X_n$
in $\C$, we have:
$$[\Delta_i P]_{X_1, \dots, X_i, X'_i, X_{i+1}, \dots, X_n}=
[P]_{X_1, \dots, X_i  \otimes X'_i, X_{i+1}, \dots, X_n}\,.$$ Now
let us prove the lemma, and assume for instance $k <i$. One may
represent $P$ and $P'$ as
$$P= \Delta_i A \Delta_i s_{k,i} B, \, P'=\Delta_i AB\,,$$
with $A \in \PB_{n-1}$, $B\in \PB_n$, and $s_{k,i}$ is the Burau
generator. Using the above-mentioned fact, we may assume $A$ and
$B$ trivial. The lemma then results from elementary properties of
the twine. The case $k> i+1$ can be treated in a similar way.
\end{PROOF}

From the lemma, we see not only that $\rvec{D}_{X_1, \dots, X_n}$
is independent of the choices made, but also that it is invariant
under Reidemeister moves of type 2 and 3. In addition, it is
invariant under `right-handed moves of type 0', namely
$$\draw{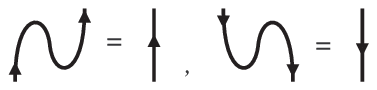}\,,$$
In the first case, it is an easy consequence of the identity
$$c_{i+1} [\Delta^3_i P]= [P]\,,$$ where $P$ is a pure braid diagram and
$\Delta^3_i P$ is is obtained by tripling  the $i$-th strand of
$P$. The second case is deduced from the first, using type 2
moves.

Now let $D$ be a arbitrary $n$-string link diagram. For each local
extremum pointing to the left, modify $D$ in the following way :
$$(1) \quad\draw{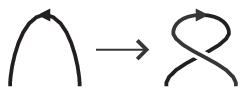}\,,\quad (2)\quad\draw{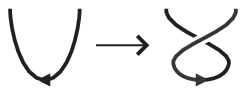}\,.$$
This operation yields a right-handed diagram $D^r$.

For $1 \le i \le n$, let $t_i$ be the algebraic number of
modifications made on the $i$-th component,  with (1) counting as
$-1$ and (2) as $+1$.

Set $\rvec{D}= (\theta^{t_1} \otimes \dots \otimes \theta^{t_n})
D^r$.

Clearly this is invariant under Reidemeister moves of type 2 and
3. As for invariance under type $0$ moves, the case when the
extrema point to the left reduces to the right-handed case
(already proved) via:
$$\draw{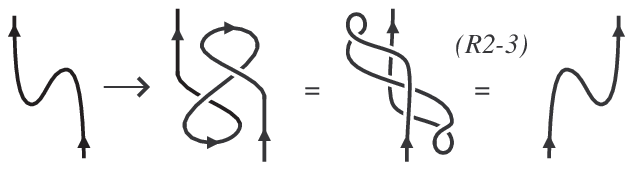}\,.$$
Moreover, we have $\raisebox{2mm}{\draw{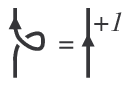}}$ (where the $+1$
denotes the twist), and $\raisebox{2mm}{\draw{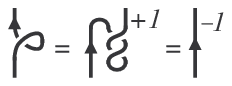}}$ by
self-duality of the twist, hence invariance under moves of type 1.

Let us summarize: given $P \in \RStL_n$, we have constructed
$\rvec{P} \in \End(\otimes^n)$.  Now notice that when one forms
the $i$-th contraction $c_i D$ of a string link diagram $D$, the
orientation of its $i$-th component is reversed; in particular, if
$D$ is right-handed,  $c_i D$ is not, unless there are no local
extrema on the $i$-th component, that is, the $i$-th component is
unknotted. In that case, we do have $c_i \rvec{P}=\rvec{c_i P}$.
Indeed, we may represent $P$ by a diagram whose $i$-th component
has no local extrema, hence $P= c_{j_m}\dots c_{j_1} Q$, with $Q$
ribbon pure braid, $j_1\le \dots \le j_m$, and $j_k \neq i+1$ for
all $k$.

\begin{LEM} The contraction operators $c_i$ satisfy the following
relations :

(a) for $i \le j-2$, $c_ic_j= c_{j-2}c_i$;

(b) for $i \ge j$, $c_ic_j= c_j c_{i+2}$.

\end{LEM}

Assume $j_{k-1} \le  i \le j_k-2$. By the lemma, $c_i P= c_{j_m-2}
\dots c_{j_k-2} c_i c_{j_{k-1}} \dots c_1 Q$, so $\rvec{c_i
P}=c_{j_m-2} \dots c_{j_k-2} c_i c_{j_{k-1}} \dots c_1 [Q]= c_i
c_{j_m}\dots c_{j_1} [Q]=c_i \rvec{P}$, hence the theorem.

\section{Turban categories}

By virtue of Shum's theorem, the category of ribbon oriented
tangles is the universal tortile category. On the other hand, we
have just seen that any category with right duals and a self-dual
twist defines invariants of ribbon string links. Recall
proposition 1:

\begin{PROP1} The category $\Turb$ (resp. $\eTang$)
is the smallest monoidal subcategory of $\Tang$ having the same
objects as $\Tang$, and containing all evaluations, coevaluations
and twists (resp. double braidings).\end{PROP1}

\medskip This suggests strongly that Shum's theorem has an analogue
for turbans. In other words, one should be able to define a notion
of `turban category', in such a way that $\Turb$ is the universal
turban category. Before we proceed to do so, let us prove
proposition 1.

\subsection*{Proof of Proposition 1}

We denote $\E$ the monoidal subcategory of $\Tang$ generated by
the evaluation morphisms. A tangle in $\E$ may be represented by a
diagram with $2n+k$ input and $k$ output, without crossings and
local min. Here is a typical example:
$$\draw{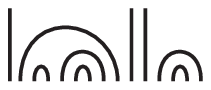}\,.$$
Dually, we denote $\E^*$ the monoidal subcategory of $\Tang$
generated by the coevaluation morphisms.

\begin{LEM} Any turban tangle $T$ may be factorized as $T=EPH$, where $P$ is a ribbon pure braid
and $E \in \E$, $H \in \E^*$. Moreover, if $T$ is even, we may
assume that each component of $P$ has trivial self-linking number.
\end{LEM}

\begin{PROOF}
Let $T$ be an (oriented) turban tangle, that is, an oriented
ribbon tangle whose linking matrix has only even entries outside
the diagonal.

We may write $T$ as $$T= \draw{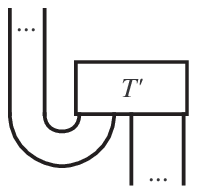}\,,$$ where $T'$ is a turban
with $2n$ input and no output. Assume $T'$ has $k$ closed
component. Pulling a local min down to the bottom line on the
right-hand side, we may represent $T'$ as
$$T'= \draw{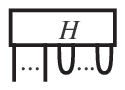}\,,$$
where $H$ is a turban with $2N=2n+2k$ input, no output and no
closed components.

Now the turban condition on $H$ excludes a configuration of four
legs $i< j< k<l$ with $(i,j)$ and $(k,l)$ connected in $H$.

By pulling one local max per component to the top line, we may
therefore write $H$ as
$$H= EQ\,,$$
where $Q \in \RStL_{2N}$, and $E$ is an element of $\E$  with $2N$
input and $0$ output.

So we may write $P=E' Q H'$, with $Q \in \RStL$, $E' \in \E$, $H'
\in \E^*$. Now any ribbon string link may be obtained from a
ribbon pure braid by a finite number of contractions, so we may
write $Q=E''PH''$, with $P \in \RPB$, $E'' \in \E$ and $H''\in
\E^*$. Setting $E=E'E''$, $H=H''H'$, we have $T=EPH$.

If $T$ is even, $Q$ be assumed even. Now by self-duality of the
twist we have $\draw{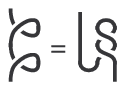}$ so we may factor $Q$ as $E'PH'$,
with $E' \in \E$, $H'\in \E^*$ and $P$ a ribbon pure braid with
trivial self-linking numbers.
\end{PROOF}

Clearly, the lemma implies proposition 1. Indeed, we know that the
category of ribbon pure braids is generated (as a monoidal
category) by the twist, and the subcategory of ribbon pure braids
with trivial self-linking numbers is generated by the double
braiding.

\medskip
\begin{REM} Proposition 1 has the following straightforward
consequences.

1) \emph{Any invariant of ribbon links or turban tangle arising
from a tortile category $\C$ is independent of the braiding:}  it
depends only on the twist and the duality.

2) If $\C$ is a twisted category with chosen right and left duals,
and $\Ob\C=\Lambda$, there exists at most one twisted,
dual-preserving functor $F_\C: \Turb[\Lambda] \to \C$ sending
$\plus_X$ ($X \in \Ob\C$) to the object $X$ itself.
\end{REM}

\medskip
Naturally, one could define a turban category to be a twisted
category with chosen left and right duals, and such that the
functor $F_\C$ exists. We would have the universal property for
free! However, such a definition would not be of great practical
use: we need a more concrete criterium. Also, it seems reasonable
to assume that the choices of left and right duals define a
sovereign structure on $\C$ (indeed, such is the case in $\Turb$).

\medskip
Let $\C$ be a sovereign category with ambidextrous twist. For any
ribbon stringlink $P$, denote $[P]=\rvec{P}=\lvec{P}$.
Pictorially, we will represent $[P]$ as $\framebox{P}\,$.
\medskip
We say that the \emph{strong sphericity condition} is satisfied if for any $P
\in \RStL_{n+2}$, $1\le i \le n+1$, $X_1, \dots, X_n, Y \in
\Ob\C$, and $f\in \End(Y)$, we have
\begin{equation}\tag{Sph}\psfrag{a}{\tiny $X_1$}
\psfrag{b}{\tiny $X_{i-1}$}\psfrag{d}{\tiny $X_i$}\psfrag{e}{\tiny
$X_n$} \psfrag{i}[c]{\tiny $Y$}\psfrag{c}{$e$}
\psfrag{f}{$\eta$}\psfrag{g}{$\eps$}\psfrag{h}{$h$}\psfrag{j}{\tiny $f$}\Draw{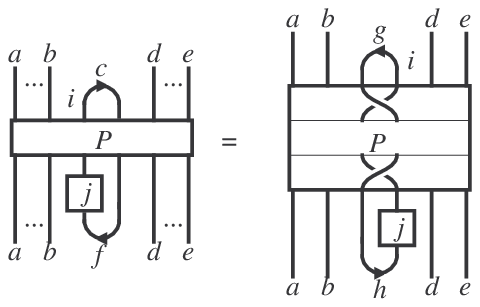}\,,\end{equation}
where $e$, $h$, $\eps$, $\eta$ denote the evaluation and
coevaluation morphisms for $Y$.
\medskip
We say that the \emph{strong interchange condition} is satisfied if for
any $P \in \RStL_{n+2}$, $1\le i \le n$, $X_1, \dots, X_n, Y \in
\Ob\C$, and $f \in \End(Y)$, we have
\begin{equation}\tag{Int}\psfrag{a}{\tiny $X_1$}
\psfrag{b}{\tiny $X_{i-1}$}\psfrag{d}{\tiny $X_i$}\psfrag{e}{\tiny
$X_{i+1}$}\psfrag{x}{\tiny $X_n$} \psfrag{i}[c]{\tiny
$Y$}\psfrag{c}{$e$}
\psfrag{f}{$\eta$}\psfrag{g}{$\eps$}\psfrag{h}{$h$}\psfrag{j}{\tiny $f$}\Draw{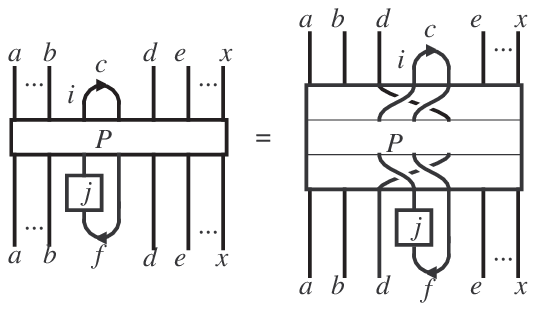}\,.\end{equation}

We say that the \emph{weak sphericity condition} (resp. the
\emph{weak interchange condition}) holds when we have (Sph) (resp. (Int))
whenever $f=\id_Y$.

\medskip

\begin{DEF}
A turban category is a twisted sovereign category with
ambidextrous twist satisfying the strong sphericity and the strong interchange
conditions.
\end{DEF}

\begin{EX}
1) For any set $\Lambda$, $\Turb[\Lambda]$ is a turban category.

2) Any tortile category is a turban category.

3) If $\C$ is a turban category, and $\D \subset \C$ is a twisted
sovereign subcategory of $\C$, then $\D$ is a turban category.
\end{EX}

\medskip
We can now state the analogue of Shum's theorem.

\begin{THEO}
The category of oriented turban tangles is the universal turban
category. In other words, if $\C$ is a turban category and
$\Lambda=\Ob\C$, there exists a unique turban functor
$$F_\C : \Turb[\Lambda] \to \C$$
sending $\plus_X$ ($X \in \Ob\C)$ to the object $X$ itself.
\end{THEO}

\begin{PROOF}

We must construct a twisted, dual-preserving functor
$$F: \Turb[\Lambda] \to \C$$ sending $\plus_X$ ($X \in \Ob\C$) to
the object $X$ itself. The proof of proposition 1 gives us a
construction for $F$, and we have to check that it is unambiguous.
We will now outline the proof.

1) The assumption that the twist is ambidextrous tells us that the
functor $F$ is well-defined on ribbon string links. It is also
well-defined on $\E$ and $\E^*$.

2) We check that $F$ is well-defined on a turban tangle $T$ with
$2n$ input, no output, and no closed component. Such a $T$ may be
factorized as $T=EP$, with $P \in \RStL_{2n}$ and $E \in \E$, so
we should set $F(T) = F(E)F(P)$. We have to check that this is
independent of the actual factorization. The proof of this fact is
similar to that of Theorem 3: starting from a suitable diagram
representing $T$, a factorization is obtained by pulling certain
local max to the top line in the right order. Just as in the proof
of theorem 3, each of these local max may have to be modified so
as to point in the appropriate direction. We need an analogue of
lemma 2, graphically:
$$\draw{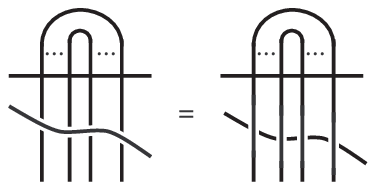}\,;$$
this is easy to check using theorem 3. This tells us that
$F(E)F(P)$ is independent of the choices of pathes, and it is then
easy to check that it depends only on the isotopy class of the
tangle.

3) We now define $F$ on turban tangles $T$ with $2n$ input, no
output, and closed components $L_1, \dots L_k$. Such a tangle may
be factorized as $$T= \draw{Handle}\,,$$ where $H$ is a turban
with $2N=2n+2k$ input, no output and no closed components. Such a
factorization is obtained by pulling a local min of each of the
$L_i$ to the bottom line, and to the right. This defines $F(T)$
with possibly two types of ambiguities: we use an numbering of
closed components, and for each closed component we must decide
whether the local minimum points to the right or to the left.
However, the (weak) exchange and sphericity conditions say precisely that
the value for $F(T)$ is independent of the numbering of
components, and the direction of each min.

4) We may write an arbitrary turban tangle $T$ as $$T=
\draw{legsup}\,,$$ where $T'$ is a turban with $2n$ input and no
output, and this defines $F$ on $T$ in an unambiguous way.

This defines a monoidal functor which has the required properties.
\end{PROOF}

\begin{REM} Theorem 4 remains true if we replace the strong sphericity and
strong interchange condition by their weak counterparts. However the strong version
will probably prove more useful. \end{REM}

\section{Construction of twisted categories}

\subsection{Toy example: the group-like case}

Let $G$ be a group and $\C=G-\vect$ the category of $G$-graded
vector spaces over a field $k$. Denote $k_g$ the simple object
consisting of one copy of $k$ in degree $g$. The dual of $k_g$ is
$k_{g\m}$, and the canonical evaluation and coevaluation morphisms
define a sovereign structure on $\C$. Each simple object has left
and right dimension equal to $1$.

An automorphism of $\otimes$ is characterized by its values on
simple objects, that is, a map $\delta :G \otimes G \to k^*$. It
is a twine if and only if $\delta(e,e)=1$ and
$\delta(g,h)\delta(gh,k)=\delta(h,k)\delta(g,hg)$ (in other words,
$\delta$ is a $2$-cocycle).

Notice that if $G$ is not commutative, $\C$ is not braided; and if
$G$ is commutative, a braiding on $\C$ corresponds to a
bicharacter $c :G \times G \to k^*$. The double braiding,
$c(h,g)c(g,h)$, is a symmetric bicharacter. Twines are far more
numerous than double braidings.

An automorphism of $\id_\C$ is given by a map $\theta : G \to
k^*$. It is a twist if and only if $\theta(e)=1$. Self-duality is equivalent
to $\theta(g\m)=\theta(g)$. Now any self-dual twist
actually defines a turban structure on $\C$.

\medskip
The invariants of ribbon links and turban tangles associated with
such turban categories contain no more information than the
linking matrix.

\subsection{Toy example: the infinitesimal case} Let $\C$ be a $k$-monoidal category,
where $k$ is a field, and define $\C[\eps]= \C\otimes_k k[\eps]$
by extending the scalars of $\C$ to the ring of dual numbers
$k[\eps]=k[X]/(X^2)$.

\medskip
Let $d_{X,Y} : X \otimes Y \to X \otimes Y$ be a functorial
morphism, $X,Y \in \C$. Set $D_{X,Y}=\id_{X\otimes Y}+\eps
d_{X,Y}$. Then $D$ is a twine on $\C[\eps]$ if and only if $d$
satisfies the following conditions :

(a) $d_{I,I}=0$;

(b) $d_{X,Y}+ d_{X\otimes Y,Z}= d_{Y,Z}+ d_{X, Y \otimes Z}$.

Let $t_X : X \to X$ be a functorial morphism, $X \in \Ob\C$. Set
$\theta_X= \id_X + \eps t_X$. Then $\theta$ is a twist if and only
if $t$ satisfies the condition $t_I=0$. If $\C$ has duals,
$\theta$ is self-dual if and only if $t_{\rdual{X}}=-t_X$.

Infinitesimal twists are expected to define turban invariants of
finite type (in this case, of degree $\le 1$).

\subsection{Tannaka theory for twined and twisted categories}

Let $k$ be a field, and $H$ a bialgebra over $k$, with coproduct
$\Delta$, counit $\eps$, product $\mu$ and unit $\eta$. Denote
$\comod H$ the monoidal category of finite dimensional right
$H$-comodules.

\medskip
\begin{DEF} A \emph{cotwinor} of $H$ is a linear form
$d : H \otimes H \to K$ satisfying the following axioms :

\begin{axioms}{codt-1}
\item[codt-1] $d$ is invertible (for the convolution product on $H
\otimes H$), and $d * \mu = \mu * d$ in $\Hom(H^{\otimes 2}, H)$;

\item[codt0] $d(\eta \otimes \eta)= 1$;

\item[codt1] $(d \otimes \eps) *  d(\mu \otimes \id_H)= (\eps
\otimes d)
* d(\id_H \otimes \mu)$;

\item[codt2] $d_{1,3}
d\m_{3,4}d_{2,4}d_{3,4}=d\m_{3,4}d_{2,4}d_{3,4}d_{1,3}$ in
$\Hom(H^{\otimes 4}, k)$.
\end{axioms}

\end{DEF}

\medskip
\begin{DEF} A \emph{cotwistor} of $H$ is a linear form $\theta : H
\to k$ satisfying the following axioms:

\begin{axioms}{cotw-1}
\item[cotw-1] $\theta$ is invertible (for the convolution product
on $H$), and $(\theta \otimes \id_H)\Delta=(\id_H \otimes
\,\theta)\Delta$;

\item[cotw0] $\theta\eta=1$;

\item[cotw1] $(\theta\m \mu)_{12}(\theta\mu^3)_{123}(\theta\m
\mu)_{23}(\theta\mu^3)_{234}(\theta\m \mu)_{34}=$

\hspace{2cm}$(\theta\m \mu)_{34}(\theta\mu^3)_{234}(\theta\m
\mu)_{23}(\theta\mu^3)_{123}(\theta\m \mu)_{12}$.
\end{axioms}

\end{DEF}

\begin{THEO}
The set of twines (resp twists) of $\comod H$ is in 1---1
correspondence with the set of cotwinors (resp. cotwistors) on
$H$. Moreover, when $H$ admits an antipode $S$ (that is, when
$\comod H$ has right duals) self-dual twists of $\comod H$
correspond exactly with cotwistors $\theta$ such that $\theta
S=\theta$.
\end{THEO}

\begin{PROOF} This is straightforward tannakan translation. Given $d: H \otimes H \to k$,
and two $H$-comodules $V$, $V'$, with coactions $\del$, $\del'$,
define $D_{V,V'}=(\id_{V \otimes V'} \otimes \, d)(\id_X \otimes
\sigma_{H,X'} \otimes \id_H)(\del \otimes \del')$. Axiom (codt-1)
means that $D_{V,V'}$ is an isomorphism of comodules, and
(codt0)-(codt2) translate axioms (DT0)-(DT2) of twines. Similarly,
given $\theta : H \to k$, and a $H$-comodule $V$, define
$\theta_V= (\id_H \otimes \theta)\del$. Axiom (cotw-1) means that
$\theta_V$ is an isomorphism of comodules, and (cotw0), (cotw1)
translate axioms (TW0), (TW1) of twists.
\end{PROOF}

Should the reader prefers modules to comodules, here are the dual notions.

\medskip
\begin{DEF} A \emph{twinor} of $H$ is an element $d \in H \otimes
H$ satisfying the following axioms :

\begin{axioms}{dt-1}
\item[dt-1] $d$ is invertible, and $\forall x \in H,
d\Delta(x)=\Delta(x)d$;

\item[dt0] $(\eps\otimes \eps)d= 1$;

\item[dt1] $(d \otimes \eta)  (\Delta \otimes \id_H)d= (\eta
\otimes d) (\id_H \otimes \Delta)d$;

\item[dt2] $d_{1,3}
d\m_{3,4}d_{2,4}d_{3,4}=d\m_{3,4}d_{2,4}d_{3,4}d_{1,3}$ in
$H^{\otimes 4}$.
\end{axioms}
\end{DEF}

\begin{DEF} A \emph{twistor} of $H$ is an element $\theta \in H$ satisfying the following axioms:

\begin{axioms}{tw-1}

\item[tw-1] $\theta$ is invertible and central;

\item[tw0] $\eps\theta=1$;

\item[tw1]
$(\Delta\theta\m)_{12}(\Delta^3\theta)_{123}(\Delta\theta\m)_{23}(\Delta^3\theta)_{234}(\Delta\theta\m)_{34}=$

\hspace{2cm}$(\Delta\theta\m)_{34}(\Delta^3\theta)_{234}(\Delta\theta\m)_{23}(\Delta^3\theta)_{123}(\Delta\theta\m)_{12}$.
\end{axioms}
\end{DEF}

\medskip
If $d$ is a twinor (resp. twistor) of $H$, the monoidal category
$H$-$\mod$ of finite dimensional left $H$-modules is entwined
(resp. twisted).

%\cite{Ver1} et \cite[page 154]{Mar1}

\bibliographystyle{amsalpha}
\bibliography{mabiblio}

\vspace{2cm} {\parindent=0pt \sc \small Alain Bruguières,

 Institut de Mathématiques et
Modélisation de Montpellier (I3M) - UMR C.N.R.S. 5149

Département des Sciences Mathématiques, Université Montpellier II,
Case Courrier 051, Place Eugène Bataillon, 34095 Montpellier Cedex
5, France

e-mail adress: \rm{bruguier@math.univ-montp2.fr}}
\end{document}